\documentclass[12pt]{amsart}
\usepackage{amssymb,enumitem,amsthm,amsmath,bibentry,thmtools,ulem,mathtools,xurl,microtype,amsaddr,tabularx}
\usepackage[pdftex,colorlinks=true,linkcolor=Black]{hyperref}
\usepackage[usenames,dvipsnames]{xcolor}

\numberwithin{equation}{section}

\theoremstyle{plain}
\newtheorem{theorem}{Theorem}[section]
\newtheorem{lemma}[theorem]{Lemma}
\newtheorem{corollary}[theorem]{Corollary}

\theoremstyle{definition}

\theoremstyle{remark}
\newtheorem{remark}[theorem]{Remark}

\numberwithin{equation}{section}

\newcommand{\ddt}{{\frac{d}{dt}}}

\newcommand{\be}{\begin{equation}}
\newcommand{\ee}{\end{equation}}
\newcommand{\bc}{\begin{dcases}}
\newcommand{\ec}{\end{dcases}}
\newcommand{\ba}{\begin{aligned}}
\newcommand{\ea}{\end{aligned} }

\title[Axisymmetric NS with total-slip boundary]{Global regularity of axisymmetric Navier--Stokes equations under total slip and critical smallness}

\author[T.L. Chan]{Tsz-Lik Chan }

\address[]{Department of Mathematics, University of California, Riverside, CA 92521,
USA}

\email{tchan204@ucr.edu}



\begin{document}

\maketitle
\begin{abstract}
We study the axisymmetric incompressible Navier--Stokes equations in a finite cylinder subject to the Navier--total--slip boundary condition, that is, the Navier slip condition with zero friction. We prove that the solution remains smooth for all time provided that the initial swirl $\Gamma_0 = r v_{\theta,0}$ has zero mean and an explicit smallness condition on the initial data holds. The condition involves $\Gamma = rv_\theta$ and the variables $V = v_\theta/\sqrt{r}$ and $\Omega = \omega_\theta/r$, and it is invariant under the natural scaling of the Navier--Stokes equations, hence critical. The total-slip condition produces a boundary term on the lateral wall that cannot be controlled by the dissipation alone; we overcome this with a sharp trace-type identity combined with a remainder controlled by $\|\Gamma\|_{L^2}$. In addition, the Robin-type boundary condition $\partial_r \Gamma = 2\Gamma$ destroys the maximum principle for $\Gamma$, and we replace it with a uniform $L^4$ bound. These estimates yield uniform $L^4$ bounds for the velocity, and regularity follows from the Ladyzhenskaya--Prodi--Serrin criterion together with boundary regularity results for the Navier slip condition.
\end{abstract}

\section{Introduction}
\subsection{Background}
The incompressible Navier--Stokes equations in $\mathbb{R}^3 \times (0,\infty)$ are
\begin{equation}
\begin{dcases}
    \partial_t v + (v \cdot \nabla) v = \nu \Delta v - \nabla P + f, \\
    \nabla \cdot v = 0,
\end{dcases}
\label{eq:NS}
\end{equation}
which describe the motion of a viscous incompressible fluid with velocity $v$, pressure $P$, external force $f$ and kinematic viscosity $\nu>0$. The first equation expresses conservation of momentum and the second conservation of mass. The viscosity is inversely proportional to the Reynolds number: at high Reynolds number the nonlinear transport term dominates, while at low Reynolds number the flow is dominated by viscous diffusion. Throughout this paper we take $\nu = 1$ and $f = 0$.

In 2D, the global regularity problem is solved by Leray \cite{leray1933etude} in $\mathbb{R}^2$, and Ladyzhenskaya \cite{ladyzhenskaya1959solution} proved that on a bounded smooth domain $\mathcal{O}\subset\mathbb{R}^2$ with the no-slip condition $v=0$ on $\partial\mathcal{O}$, every divergence-free $v_0\in H^1_0(\mathcal{O})$ and $f\in L^2(0,T;L^2(\mathcal{O}))$, $\partial_t f \in L^2(\mathcal{O})$, give a unique global regular solution. 

In 3D, Leray \cite{leray1934mouvement} constructed global weak solutions on $\mathbb{R}^3$, and Hopf \cite{hopf1950anfangswertaufgabe} extended the construction to bounded domains. These Leray--Hopf weak solutions satisfy $v\in L^\infty(0,T;L^2_\sigma)\cap L^2(0,T;H^1)$ for all $T>0$ and the energy inequality
\begin{equation}
E(v)(t) := \|v(t)\|_{L^2}^2 + 2\nu \int_0^t \|\nabla v(s)\|_{L^2}^2 \, ds \le \|v_0\|_{L^2}^2 .
\end{equation}
Whether they remain smooth for all time, or whether singularities can form in finite time from smooth data, is a major open problem.

The difficulty is reflected in the scaling of the equations. If $v$ solves \eqref{eq:NS} with $f=0$, then so does
\begin{equation}
v_\lambda(x,t) = \lambda v(\lambda x, \lambda^2 t), \qquad P_\lambda(x,t) = \lambda^2 P(\lambda x, \lambda^2 t),
\end{equation}
and in $\mathbb{R}^n$ the energy transforms as $E(v_\lambda) = \lambda^{2-n} E(v)$. For $n=2$ the energy is scale invariant (critical). For $n=3$ we have $E(v_\lambda)=\lambda^{-1}E(v)$. Zooming into a region corresponds to sending $\lambda\to 0$, and then $E(v_\lambda)\to\infty$: at small scales the rescaled solution has large energy, so the energy bound gives no control of the small-scale behavior, and the problem is supercritical. 

Regularity criteria are therefore usually stated in terms of scale-invariant quantities. The classical example is the Ladyzhenskaya--Prodi--Serrin criterion \cite{serrin1962interior,prodi1959teorema,ladyzhenskaya1967uniqueness}: if a weak solution satisfies
\begin{equation}
v \in L^{p}(0,T; L^{q}), \qquad \frac{2}{p} + \frac{3}{q} \le 1, \quad 3 < q \leq \infty,
\end{equation}
then $v$ is smooth on $(0,T]$; the endpoint $L^\infty(0,T;L^3)$ was established in \cite{escauriaza2003}.

A natural intermediate setting is the axisymmetric flows, which are invariant under rotations about a fixed axis; see for example \cite{chae2002regularity,chen2009lower1,chen2009lower2,hou2008global,lei2017criticality,liu2009characterization,zhang2022bounded}. With $r=\sqrt{x_1^2+x_2^2}$ and the cylindrical basis
\begin{equation}
e_r = \Big(\frac{x_1}{r}, \frac{x_2}{r}, 0\Big), \qquad e_\theta = \Big( -\frac{x_2}{r}, \frac{x_1}{r}, 0\Big), \qquad e_3 = (0,0,1),
\end{equation}
an axisymmetric velocity has the form
\begin{equation}
v = v_r(r,x_3,t)\,e_r + v_\theta(r,x_3,t)\,e_\theta + v_3(r,x_3,t)\,e_3 .
\end{equation}
Historically, it has been known since the work of Ladyzhenskaya, Ukhovskii, and Iudovich that in the \textit{swirl-free} case \(v_\theta \equiv 0\), finite-energy solutions are globally smooth \cite{ukhovskii1968axially,ladyzhenskaya1968unique}. With swirl it remains open, but the swirl $rv_\theta$ satisfies a transport--diffusion equation and is scale invariant. This has led to regularity results under smallness of $rv_\theta$; see \cite{lei2017criticality}. 

\subsection{Problem setting and main result}
In this paper we study axisymmetric solutions of \eqref{eq:NS} with $\nu=1$ and $f=0$ in the finite cylinder
\be
D=\{ (r,\theta,x_3): 0\leq r \leq 1,\ 0 \leq \theta < 2\pi,\ 0 \leq x_3 \leq 1 \}.
\ee
From now on we denote the velocity by $v = v_r e_r + v_\theta e_\theta + v_3 e_3$ and we write $b = v_r e_r + v_3 e_3$. We write the boundary of $D$ as the union of the horizontal boundary $\partial^H D$ and the vertical boundary $\partial^V D$. We suppose that $v$ satisfies the Navier--total--slip boundary condition 
\be
v\cdot n = 0, \quad [\mathbb{S}(v)\, n]_{tan} = 0,
\ee
where $\mathbb{S}(v)=\frac{1}{2}(\nabla v + (\nabla v)^T)$ is the symmetric part of $\nabla v$ and $[\mathbb{S}(v)\, n]_{tan}$ means the tangential component of $\mathbb{S}(v)\, n$. 

We want to mention that the Navier--total--slip boundary condition has been studied  before. For example, the recent paper \cite{li2024finite} studied this boundary condition in the exterior of a cone and also \cite{li2025pipe} studied this boundary condition in a infinitely long pipe.

The Navier--total--slip boundary condition in cylindrical coordinates reads
\be \ba
& v_3 = 0, \quad \partial_{x_3}v_r = \partial_{x_3} v_\theta = 0, \quad \text{on $\partial^H D \times (0,T],$} \\
& v_r = 0, \quad \partial_{r}v_\theta = \frac{v_\theta}{r} , \quad \partial_{r} v_3 = 0,\quad \text{on $\partial^V D \times (0,T].$}
\label{eq:U_BC}
\ea \ee
Writing out \eqref{eq:NS} in the cylindrical coordinate basis, we find that $v$ satisfies
\be \bc
    \label{eq:ASNS}
	&(\Delta - \frac{1}{r^2}) v_r - (b\cdot\nabla) v_r + \frac{v_\theta^2}{r} - \partial_r P - \partial_t v_r = 0, \\
	&(\Delta - \frac{1}{r^2}) v_\theta - (b\cdot\nabla) v_\theta - \frac{v_\theta v_r}{r} - \partial_t v_\theta =  0, \\
	&\Delta v_3 - (b\cdot\nabla) v_3 - \partial_{x_3} P - \partial_t v_3 = 0, \\
	& \frac{1}{r} \partial_r (r v_r) + \partial_{x_3}v_3 =0, \\
	&\partial_{x_3}v_r = \partial_{x_3} v_\theta = 0, \quad v_3 =0 \quad \text{on $\partial^H D \times (0,T],$} \\
	&\partial_{r}v_\theta = \frac{v_\theta}{r} , \quad \partial_{r} v_3 = 0,\quad v_r = 0 \quad \text{on $\partial^V D \times (0,T],$}
\ec \ee
where $\Delta = \partial_r^2 + \frac{1}{r}\partial_r + \partial_{x_3}^2$, and the initial condition
\begin{equation}
v(x,0) = v_0(x), \quad x \in D.
\end{equation}

Let $\omega = \operatorname{curl} \, v $ be the vorticity. In cylindrical coordinates it is given by
\begin{equation}
 \omega = (-\partial_{x_3} v_\theta) e_r + (\partial_{x_3} v_r - \partial_r v_3) e_\theta + \Big( \frac{v_\theta}{r} + \partial_r v_\theta\Big) e_3,
\end{equation}
and the vorticity equation $\Delta \omega - (v \cdot \nabla) \omega + ( \omega \cdot  \nabla ) v - \partial_t \omega = 0$ becomes
\be \bc
    \label{eq:vorticity}
	&(\Delta - \frac{1}{r^2}) \omega_r - (b\cdot\nabla) \omega_r + (\omega_r \partial_r + \omega_3 \partial_{x_3} )v_r  - \partial_t \omega_r = 0, \\
	&(\Delta - \frac{1}{r^2}) \omega_\theta - (b\cdot\nabla) \omega_\theta + \frac{2}{r} v_\theta \partial_{x_3}v_\theta  + \frac{\omega_\theta v_r}{r} - \partial_t \omega_\theta =  0, \\
	&\Delta \omega_3 - (b\cdot\nabla) \omega_3 +(\omega_r \partial_r + \omega_3 \partial_{x_3}) v_3  - \partial_t \omega_3 = 0.
\ec \ee

The swirl $\Gamma = rv_\theta$ satisfies
\begin{equation}
	\Delta \Gamma - (b\cdot\nabla) \Gamma - \frac{2}{r} \partial_r \Gamma - \partial_t \Gamma = 0,
    \label{eq:GammaU}
\end{equation}
with the boundary conditions
\be \bc
\partial_{x_3}\Gamma = 0 \quad \text{on $\partial^H D$,} \\
\partial_r \Gamma = 2\Gamma \quad \text{on $\partial^V D$,}
\label{eq:GammaU_bc}
\ec \ee
where the second condition follows from $\partial_r \Gamma = v_\theta + r\partial_r v_\theta = 2 v_\theta$ at $r=1$. From the equation of $v_\theta$ in \eqref{eq:ASNS} and the equation of $\omega_\theta$ in \eqref{eq:vorticity}, $V = v_\theta/\sqrt{r}$ and $\Omega=\omega_\theta/r$ satisfy
\be \ba
\label{eq:VOmega}
&\Delta V + \frac{1}{r}\partial_r V - b\cdot \nabla V - \frac{3}{4r^2}V  - \frac{3v_r}{2r} V - \partial_t V = 0,\\
&\Delta \Omega + \frac{2}{r} \partial_r \Omega - b\cdot\nabla \Omega + \frac{1}{r} \partial_{x_3} (V^2) - \partial_t \Omega = 0,
\ea \ee
with the boundary conditions
\be \ba
\label{eq:VOmega_bc}
& \Omega = 0 \quad \text{on $\partial D$,} \\
& \partial_{x_3} V = 0 \quad \text{on $\partial^H D$,} \\
& \partial_r V = \frac{1}{2r} V \quad \text{on $\partial^V D$.}
\ea \ee

We point out that apart from the $V$, $\Omega$ system, one can also formulate a closed system in terms of $\Omega$ and $J = \omega_r/r$ for proving global regularity results \cite{lei2017criticality,wei2016regularity}. We follow the method shown in \cite{lei2017criticality} to formulate a closed system in terms of $\Omega$ and $V$ because of its simplicity relative to the formulation in terms of $\Omega$ and $J = \omega_r/r$. 

We also recall from \cite{zhang2022bounded} if $v$ is an axisymmetric vector then the formula for $\nabla v$ is:
\begin{equation}
\nabla v = 
\begin{pmatrix}
    \partial_r v_r & -\frac{1}{r} v_\theta & \partial_{x_3} v_r \\
\partial_r v_\theta & \frac{1}{r} v_r & \partial_{x_3} v_\theta \\
\partial_r v_3 & 0 & \partial_{x_3} v_3
\end{pmatrix}.
\end{equation}

\subsection{Notation}\label{sec:notation}
Here we briefly describe some notation. We will denote $b = v_r e_r + v_3 e_3.$ To simplify the notation, We omit the factor $2\pi$ coming from the integration in $\theta$ in all integrals over $D$, i.e. $\int_D f\,dx := \int_0^1\int_0^1 f\, r\,drdx_3$ for axisymmetric $f$. In particular we write $|D| = \int_0^1\int_0^1 r \,drdx_3 = \frac12$, whereas the actual volume of $D$ is $\pi$.

We let $\|f\|_{L^p(E)}$ denotes the $L^p$ norm for $f$ over a set $E$, i.e. $(\int_E |f|^p \,dx)^{1/p} $. If $E$ is omitted, it means that we are integrating over the set $D$. We use $f_0$ or $f_{\theta,0}$ to denote the initial value for $f$ or some component of $f$ if $f$ is a vector.

We also let
\be
    C_S = C_{sob} \Big( 1 + \frac{1}{\pi} \Big)
\ee
where $C_{sob}$ is the constant in the Sobolev inequality in $D$.

\subsection{Main theorem}
The main theorem is
\begin{theorem}
    Let $C_2$ be the constant in Lemma \ref{lemma:agmon_inequality}. Define
    \be \ba
        G &:= \Big( \max \Big\{    \|\Gamma_{0}\|_{L^4}^4,\ C_4 \| \Gamma_{0} \|_{L^2}^4  \Big\} \Big)^{1/4}, \\
        C_4 &:= \frac{4}{3\pi^2} \Big(  4\cdot \frac{87}{5}\cdot \Big(\frac{174}{5}\Big)^{3/2} C_S^3 + 87 + \frac{3}{2} \pi^2 \Big),
        \label{eq:G}
    \ea \ee
    and
    \be
        F := \max \Big\{ \frac{1}{4} \|V^2_0\|_{L^2}^2 + \frac{1}{40} \|\Omega_0\|_{L^2}^2,\ 54 \| \Gamma_{0} \|_{L^2}^4 \Big\}.
        \label{eq:F}
    \ee
    Suppose that
    \begin{enumerate}
        \item $v_0$ is smooth, divergence--free and satisfies the boundary condition \eqref{eq:U_BC},
        \item $\displaystyle \int_D \Gamma_{0} \, dx = 0$,
        \item $\displaystyle 80\cdot \frac{2187}{4096}\cdot 400^3\, C_2^4\, F\, G^4 \leq \frac{1}{160}$.
    \end{enumerate}
   Then $v(t)$ is smooth for all $t \geq 0$. Moreover, for all $t\ge 0$,
   \be
   \frac14\|V^2(t)\|_{L^2}^2 + \frac1{40}\|\Omega(t)\|_{L^2}^2 \le F .
   \ee
   \label{thm:main}
\end{theorem}

The smallness condition (3) in Theorem \ref{thm:main} is stated in terms of the auxiliary quantities $F$ and $G$, which are defined through maxima. The next corollary shows that it is implied by a simpler condition involving only two natural quantities of the initial data: the weighted energy
\be
f(0) := \frac{1}{4}\|V^2_0\|_{L^2}^2 + \frac{1}{40}\|\Omega_0\|_{L^2}^2
\ee
and the $L^4$ norm of the initial swirl $\Gamma_0$. The price is a smaller numerical constant. The condition is of product form, so $f(0)$ may be large provided $\|\Gamma_0\|_{L^4}$ is correspondingly small.

\begin{corollary}
\label{cor:simple}
Suppose that conditions (1)--(2) of Theorem \ref{thm:main} hold, and that
\be
f(0)\, \|\Gamma_0\|_{L^4}^4 \le \Big( 80^2 \cdot 108 \cdot C_4 \cdot \frac{2187}{4096}\cdot 400^3\, C_2^4 \Big)^{-1},
\label{eq:simple_condition}
\ee
where $C_2$ and $C_4$ are as in Theorem \ref{thm:main}. Then $v(t)$ is smooth for all $t \ge 0$.
\end{corollary}

\begin{proof}
We show that \eqref{eq:simple_condition} implies condition (3) of Theorem \ref{thm:main}.

\emph{Bound for $G$.} By H\"older's inequality and $|D| = \frac12$,
\be
\|\Gamma_0\|_{L^2}^4 \le |D|\, \|\Gamma_0\|_{L^4}^4 = \frac{1}{2}\|\Gamma_0\|_{L^4}^4 .
\ee
Since $C_4 \ge \frac{4}{3\pi^2}\cdot 87 > 2$, this gives
\be
G^4 = \max\Big\{ \|\Gamma_0\|_{L^4}^4,\ C_4\|\Gamma_0\|_{L^2}^4 \Big\} \le \max\Big\{1, \frac{C_4}{2}\Big\}\|\Gamma_0\|_{L^4}^4 = \frac{C_4}{2}\|\Gamma_0\|_{L^4}^4 .
\ee

\emph{Bound for $F$.} Since $v_\theta = \sqrt{r}\,V$, we have $\Gamma^4 = r^4 v_\theta^4 = r^6 V^4 \le V^4$ in $D$, so $\|\Gamma_0\|_{L^4}^4 \le \|V^2_0\|_{L^2}^2$. Together with the H\"older bound above,
\be
54\|\Gamma_0\|_{L^2}^4 \le 27 \|\Gamma_0\|_{L^4}^4 \le 27 \|V^2(0)\|_{L^2}^2 \le 108\, f(0).
\ee
Hence
\be
F = \max\Big\{ f(0),\ 54\|\Gamma_0\|_{L^2}^4 \Big\} \le 108\, f(0).
\ee

Combining the two bounds with \eqref{eq:simple_condition},
\be
80\cdot \frac{2187}{4096}\cdot 400^3\, C_2^4\, F\, G^4 \le 80\cdot 108\cdot \frac{C_4}{2}\cdot \frac{2187}{4096}\cdot 400^3\, C_2^4\, f(0)\,\|\Gamma_0\|_{L^4}^4 \le \frac{1}{160},
\ee
which is condition (3) of Theorem \ref{thm:main}. The corollary follows.
\end{proof}

\begin{remark}[Criticality]
\label{remark:scaling}
The condition \eqref{eq:simple_condition} is invariant under the natural scaling of the Navier--Stokes equations. From the natural scaling $\tilde v(x,t) = \lambda v(\lambda x, \lambda^2 t)$ we have $\tilde V(x) = \lambda^{3/2}V(\lambda x)$, $\tilde\Omega(x) = \lambda^{3}\Omega(\lambda x)$, $\tilde\Gamma(x) = \Gamma(\lambda x)$, and $\int g(\lambda x)\,dx = \lambda^{-3}\int g(x)\,dx$, we have
\be
\|\tilde V^2_0\|_{L^2}^2 = \lambda^{3}\|V^2_0\|_{L^2}^2, \qquad
\|\tilde \Omega_0\|_{L^2}^2 = \lambda^{3}\|\Omega_0\|_{L^2}^2, \qquad
\|\tilde \Gamma_0\|_{L^4}^4 = \lambda^{-3}\|\Gamma_0\|_{L^4}^4 .
\ee
Therefore $\tilde f(0) = \lambda^3 f(0)$, and the product $f(0)\,\|\Gamma_0\|_{L^4}^4$ is invariant. The zero-mean condition $\int \Gamma_0\, dx = 0$ is also invariant. Thus \eqref{eq:simple_condition} is critical. This is analogous to the criticality of the condition in \cite[Remark 1.2]{chan2026global}.
\end{remark}

\begin{remark}[On the zero-mean condition]
\label{remark:mean}
By Lemma \ref{lemma:int_Gamma_conserved} below, the angular momentum $\int_D \Gamma \, dx$ is conserved, so condition (2) holds for all time. Condition (2) is used to obtain the uniform bounds for $\|\Gamma\|_{L^2}$ and $\|\Gamma\|_{L^4}$ in Lemmas \ref{lemma:L^2_Gamma_U} and \ref{lemma:L^4_Gamma_U}, where the Robin-type condition $\partial_r\Gamma = 2\Gamma$ on $\partial^V D$ prevents a direct maximum principle. Under the total-slip condition every rigid rotation $c\, r e_\theta$ is a stationary solution, with $\Gamma = c r^2$ and $\int_D \Gamma\,dx = c/4 \neq 0$; these are excluded by condition (2).
\end{remark}

\subsection{Organization of the paper}
The rest of the paper is organized as follows. Section \ref{sec:prerequisites} collects the 
estimates needed in the proof. In Section \ref{sec:poincare} we show that the angular momentum 
$\int_D \Gamma\,dx$ is conserved and prove a Poincar\'e inequality with constant $1/\pi^2$ for 
mean-zero axisymmetric functions in $D$. Section \ref{sec:Gamma} is devoted to uniform bounds 
on $\Gamma$. Because of the Robin-type condition $\partial_r\Gamma = 2\Gamma$ on $\partial^V D$, 
the maximum principle is not available, and the energy identity for $\Gamma$ contains a boundary 
term on the lateral wall. We write this boundary term as an interior integral using the cutoff 
function $\eta(r)=\sin(\pi r/2)$. The zero-mean condition then allows the Poincar\'e inequality 
to absorb it into the dissipation, which yields $\|\Gamma(t)\|_{L^2}\le\|\Gamma_0\|_{L^2}$. 
The same argument applied to $\Gamma^2$, combined with the Sobolev inequality, gives a uniform 
$L^4$ bound for $\Gamma$ in place of the maximum principle. In Section \ref{sec:known} we recall 
estimates for $v_r/r$ from \cite{chan2026global} and check that they remain valid under the 
total-slip condition. Section \ref{sec:trace} contains the main new ingredient, a trace-type 
inequality for the boundary term $\int_0^1 V^4(1,x_3,t)\,dx_3$ that appears in the energy 
estimate for $V$. The trace inequality is proved by first noticing the identity
\[
\int_D |\nabla V^2|^2 + \frac{V^4}{r^2}\,dx - \int_0^1 V^4(1,x_3,t)\,dx_3 = \int_D r^2|\nabla \phi|^2\,dx,
\qquad \phi = \frac{V^2}{r},
\]
and then we split the boundary value of $\phi$ into its weighted mean, which equals $6\|\Gamma\|_{L^2}^2$, 
and a fluctuation controlled by $\int_D r^2|\nabla \phi|^2\,dx$.

Section \ref{sec:proof} contains the proof of Theorem \ref{thm:main}. After recalling the local existence of a strong solution, we derive a coupled energy estimate for $V$ and $\Omega$. The trace inequality controls the boundary term, the remainder is bounded by 
$\|\Gamma_0\|_{L^2}^4$, and the nonlinear term $\int_D \frac{v_r}{r}V^4\,dx$ is estimated by the Agmon-type inequality for $v_r/r$ and the $L^4$ bound for $\Gamma$. Under the smallness condition (3), a continuity argument shows that 
$\frac14\|V^2\|_{L^2}^2 + \frac1{40}\|\Omega\|_{L^2}^2$ stays below $F$ for all time. We then show $\|\nabla b\|_{L^2}\le\|\Omega\|_{L^2}$, which gives a uniform $L^4$ bound for $b = v_r e_r + v_3 e_3$. 

In Section \ref{sec:completion} we combine this with the $L^4$ bound 
for $v_\theta$ and conclude by the Ladyzhenskaya--Prodi--Serrin criterion together with the 
boundary regularity result of \cite{neustupa2018regularity}.

\section{Prerequisites}
\label{sec:prerequisites}

\subsection{Conservation law and a Poincar\'e inequality}\label{sec:poincare}
The next lemma says that the quantity $\int_D \Gamma \,dx$ is conserved.
\begin{lemma}
    If $\Gamma$ is a smooth solution to \eqref{eq:GammaU} and \eqref{eq:GammaU_bc}, then
    \be
    \ddt \int_D \Gamma \,dx = 0.
    \ee
    \label{lemma:int_Gamma_conserved}
\end{lemma} 
\begin{proof}
    We start with the equation for $\Gamma$:
    \be
    \Delta \Gamma - (b \cdot \nabla) \Gamma - \frac{2}{r}\partial_r \Gamma - \partial_t \Gamma = 0.
    \ee
    Integrating in space, we get 
    \be
    \ddt \int_D \Gamma \,dx = \int_D \Delta \Gamma \, dx - \int_D b \cdot \nabla \Gamma\, dx - 2 \int_D \frac{1}{r} \partial_r \Gamma \, dx =: I_1 + I_2 + I_3.
    \ee
    We want to show that the right side vanishes. First we can compute using integration by parts to find that
    \be \ba
    -I_2 = \int_D b \cdot \nabla \Gamma \, dx &= \int_0^1 \int_0^1 v_r \partial_r\Gamma \,r drdx_3 + \int_0^1 \int_0^1 v_3 \partial_{x_3}\Gamma \, r drdx_3 \\
    & = \int_0^1 \int_0^1 \partial_r( v_r \Gamma r) - \partial_r v_r \Gamma r - v_r \Gamma \, drdx_3 + \int_0^1\int_0^1 \partial_{x_3}(v_3 \Gamma r) - \partial_{x_3}v_3 \Gamma r \, drdx_3 \\
    & = \int_0^1 \int_0^1 \partial_r(v_r \Gamma r) + \partial_{x_3}( v_3 \Gamma r) \, drdx_3 -\int_0^1\int_0^1 \left( \partial_r v_r + \frac{v_r}{r} + \partial_{x_3}v_3 \right) \Gamma r drdx_3 \\
    & = \int_0^1 \int_0^1 \partial_r(v_r \Gamma r) + \partial_{x_3}( v_3 \Gamma r) - \int_0^1\int_0^1 \left( \nabla\cdot b \right) \Gamma r drdx_3 \\
    & = \int_0^1 v_r \Gamma r \Big|^{r=1}_{r=0} \, dx_3 + \int_0^1 v_3 \Gamma r \Big|_{x_3=0}^{x_3=1} \, dr - \int_0^1\int_0^1 \left( \nabla\cdot b \right) \Gamma r drdx_3
    \\ & = 0,
    \ea \ee
    where the last term vanishes because of the boundary conditions: $v_r=0$ on $\partial^V D$, $v_{3}=0$ on $\partial^H D$ and $\nabla \cdot b = 0$.
    Next, we compute 
    \be \ba
    I_1 := \int_D \Delta \Gamma \, dx  &= \int_0^1\int_0^1 \partial_r^2 \Gamma \, rdrdx_3 + \int_0^1\int_0^1 \frac{1}{r}\partial_r \Gamma \, rdrdx_3 + \int_0^1\int_0^1 \partial_{x_3}^2 \Gamma \, rdrdx_3 \\
    & = T_1 + T_2 + T_3.
    \ea \ee
    Using integration by parts we have 
    \be \ba
    T_1 &= \int_0^1\int_0^1 \partial_r( \partial_r \Gamma r) - \partial_r \Gamma \, drdx_3 \\
    & = \int_0^1 \partial_r \Gamma r \Big|_{r=0}^{r=1}  dx_3 - \int_0^1 \Gamma \Big|_{r=0}^{r=1} dx_3 \\
    & = \int_0^1 \partial_r \Gamma(1,x_3,t) \, dx_3 - \int_0^1 \Gamma(1,x_3,t) \, dx_3,
    \ea \ee
    \be \ba
    T_2 = \int_0^1 \int_0^1 \partial_r \Gamma \, drdx_3 = \int_0^1 \Gamma(1,x_3,t) \, dx_3,
    \ea \ee
    where $\Gamma(0,x_3,t)$ vanishes due to the regularity of $v_\theta$ \cite[Corollary 1]{liu2009characterization}.
    \be
    T_3 = \int_0^1 r\partial_{x_3} \Gamma(r,1,t) - r \partial_{x_3} \Gamma(r,0,t) \, dr
    \ee
    Due to the boundary condition $\partial_{x_3}\Gamma = r \partial_{x_3}v_\theta = 0$ on $\partial^H D$, we have $T_3 = 0$. Then
    \be
    I_1 = T_1+T_2+T_3 = \int_0^1 \partial_r \Gamma(1,x_3,t) \, dx_3 = 2\int_0^1 \Gamma(1,x_3,t) dx_3. \
    \ee
    For $I_3$ we have,
    \be
    I_3 = -2\int_0^1 \int_0^1  \partial_r \Gamma \, drdx_3 = -2 \int_0^1 \Gamma(1,x_3,t) \, dx_3.
    \ee
    Finally, we get 
    \be
    \ddt \int_D \Gamma dx = I_1 + I_2 + I_3 = 0.
    \ee
\end{proof}

From Lemma \ref{lemma:int_Gamma_conserved}, we know that if $\int_D \Gamma_{0} \,dx =0$, then $\int_D \Gamma\,dx=0$ for any $t\geq0$.

Next, let 
\be \ba
&A_m = \Big\{ Q(r) \in H^1(1/m,1) \backslash \{0 \}: \int_{1/m}^1 Q(r) \, rdr = 0 \Big\},
\\ & A = \Big\{ R(r): \int_0^1 rR^2  + r(R^\prime)^2 \, dr < \infty, R \neq 0, \int_{0}^1 R(r) \, rdr = 0 \Big\},
\ea \ee
and let 
\be
\lambda_{r,m} = \inf_{Q \in A_m } \frac{\int_{1/m}^1 (Q'(r))^2 \, rdr}{\int_{1/m}^1 (Q(r))^2 \, rdr}.
\ee
By Lemma 2.3 in \cite{li2024finite}, we know that
\be
\lambda_{r,m} \geq \frac{\pi^2}{\Big(1-\frac{1}{m}\Big)^2} - \max_{ r \in [1/m,1] } \left( \frac{1}{2} \frac{(r)''}{r} - \frac{3}{4} \frac{((r)')^2}{(r)^2} \right) \geq \pi^2 + \frac{3}{4}.
\label{eq:ray_1}
\ee
We wish to let $m \rightarrow \infty$. First suppose $R\in A$ and define
\be
c_m = \Big( \int_{1/m}^1 r\,dr \Big)^{-1} \int_{1/m}^1 R(r)\,r dr.
\ee
By the Triangle inequality,
\be
   \Big( \int_{1/m}^1 R^2(r) \, rdr \Big)^{1/2} \leq \Big( \int_{1/m}^1 (R(r)-c_m)^2 \, rdr \Big)^{1/2} + |c_m| \Big( \int_{1/m}^1 r\,dr \Big)^{1/2}.
   \label{eq:ray_2}
\ee

For the first term on the right side, note that
\be
\int_{1/m}^1 (R^2 + R'^2)\,dr \le m\int_0^1 (R^2+R'^2)\,r\,dr < \infty,
\ee
so $R\in H^1(1/m,1)$, and hence $R-c_m\in H^1(1/m,1)$ with $(R-c_m)'=R'$. By the definition of $c_m$, $\int_{1/m}^1 (R-c_m)\,r\,dr=0$, so $R-c_m\in A_m$ (if $R-c_m\equiv 0$ on $[1/m,1]$, the estimate below is trivial). Therefore \eqref{eq:ray_1} gives
\be
    \int_{1/m}^1  ( R(r) - c_m)^2 \, rdr \leq \frac{1}{\pi^2 + 3/4} \int_{1/m}^1 (R'(r))^2 \, rdr \leq \frac{1}{\pi^2 + 3/4} \int_{0}^1 (R'(r))^2 \, rdr.
\ee

Next we handle the second term on the right side of \eqref{eq:ray_2}, we want to show that $|c_m| \rightarrow 0$ as $m \rightarrow \infty$. First by the zero mean condition in $A$ we get
\be
- \int_{0}^{1/m} R(r) \, rdr = \int_{1/m}^1 R(r)\, rdr.
\ee
Then by Cauchy–Schwarz inequality
\be
    \Big|  \int_{1/m}^1 R(r)\, rdr \Big| = \Big| \int_{0}^{1/m} R(r) \, rdr \Big| \leq \Big( \int_0^{1/m} r \, dr \Big)^{1/2} \Big( \int_0^1 R^2(r)\, rdr \Big)^{1/2} \rightarrow 0 \, \text{as $m \rightarrow \infty$,}
\ee
since $\int_0^1 R^2 r\,dr < \infty.$ And since $\int_{1/m}^1 r\,dr  \rightarrow 1/2$, we have as $m \rightarrow \infty$, $|c_m|\rightarrow 0$ and
\be
|c_m| \Big( \int_{1/m}^1 r\,dr \Big)^{1/2} \rightarrow 0.
\ee
For \eqref{eq:ray_2}, the left side converges to $(\int_0^1 R^2 \, rdr)^{1/2}$ by the Monotone Convergence Theorem. By taking the limit in \eqref{eq:ray_2}, we have for any $R \in A$,
\be
    \int_{0}^1 R^2(r) \, rdr \leq \frac{1}{\pi^2 + 3/4} \int_{0}^1 (R'(r))^2 \, rdr.
   \label{eq:ray_3}
\ee

For the $x_3$ direction, let
\be
B = \Big\{ S(x_3) \in H^1(0,1) \backslash \{ 0 \} : \int_0^1 S(x_3) dx_3 = 0 \Big\},
\ee
and
\be
\lambda_3 = \inf_{S \in B} \frac{\int_0^1 (S'(x_3))^2 \,dx_3}{\int_0^1 (S(x_3))^2 \,dx_3}.
\ee
Again using Lemma 2.3 in \cite{li2024finite} we know that
\be
\lambda_3 \geq \pi^2 - \max_{ [0,1] } \left( \frac{1}{2} \frac{(1)''}{1} - \frac{3}{4} \frac{((1)')^2}{1} \right) = \pi^2.
\ee

A version of the next lemma appears in \cite[Lemma 4.3]{li2026cone} for a special domain in the spherical coordinate. We prove it again for clarity since the domain in this paper is the finite cylinder in the cylindrical coordinate.
\begin{lemma}
    \label{lemma:Poincare}
    Suppose $f(r,x_3) \in H^1(D)$ is an real-valued axisymmetric scalar function function that satisfies 
    \be
    \int_D f \, dx = 0,
    \ee
    then 
    \be
    \int_D |f|^2 \, dx \leq \frac{1}{\pi^2} \int_D |\nabla f|^2 \, dx.
    \ee
\end{lemma}
\begin{proof}
    Define 
    \be
        g(r) = \int_0^1 f(r,x_3) \, dx_3, \quad h(r,x_3) = f(r,x_3) - g(r).
    \ee
    Then $f = g + h$. 
    
    By the fact that 
    \be
        \int_0^1 \int_0^1 f\, rdrdx_3 = 0,
    \ee
    we see that
    \be
        \int_0^1 g(r) \, r dr = 0, \quad \int_0^1 h(r,x_3) \, dx_3 = \int_0^1 f(r,x_3) \, dx_3 - \int_0^1 g(r) \, dx_3 = 0.
    \ee
    Thus $g \in A$ and $h(r,\cdot) \in B$ and hence we have the following inequalities
    \be \ba
        &\int_0^1 g^2 \, rdr \leq \frac{1}{\pi^2 + 3/4} \int_0^1 |\partial_r g|^2 \, rdr \\
        &\int_0^1 h^2 dx_3 \leq  \frac{1}{\pi^2} \int_0^1 |\partial_{x_3}h|^2 \, dx_3.
    \ea \ee
    On the other hand we have 
    \be
        \int_0^1 \int_0^1  g(r) h(r,x_3) \,rdrdx_3 = \int_0^1 \int_0^1 h(r,x_3)\, dx_3  \, g(r) \, rdr = 0.
        \label{eq:lemma_Poincare_1}
    \ee
    and since for any $r$, $\int_0^1 h\, dx_3 = 0$,
    \be
        \int_0^1 \int_0^1 (\partial_r g)(\partial_r h) \,rdrdx_3 = \int_0^1 \int_0^1 \partial_r h\, dx_3 \, \partial_r g \, rdr = \int_0^1 \partial_r \int_0^1 h\, dx_3 \, \partial_r g \, rdr = 0,
    \ee
    This implies that
    \be
    \int_0^1 \int_0^1 |\partial_r f|^2 \, rdrdx_3 = \int_0^1 \int_0^1 |\partial_r g|^2 \, rdrdx_3 + \int_0^1 \int_0^1 |\partial_r h|^2 \, rdrdx_3 \geq \int_0^1 \int_0^1 |\partial_r g|^2 \, rdrdx_3.
    \label{eq:lemma_Poincare_2}
    \ee
    Also by the definition of $f,g,h$, 
    \be
    \int_0^1 \int_0^1 |\partial_{x_3} f|^2 \, rdrdx_3 = \int_0^1 \int_0^1 |\partial_{x_3} h|^2 \, rdrdx_3.
    \label{eq:lemma_Poincare_3}
    \ee
    Hence, using \eqref{eq:lemma_Poincare_1}, \eqref{eq:lemma_Poincare_2}, \eqref{eq:lemma_Poincare_3} and the fact that $f= g+h$, we obtain
    \be \ba
        \| f \|_{L^2}^2 &= \|g\|_{L^2}^2 + \|h\|_{L^2}^2  \\
        & \leq \frac{1}{\pi^2 + 3/4} \| \partial_r g \|_{L^2}^2 + \frac{1}{\pi^2} \| \partial_{x_3} h \|_{L^2}^2  \\
        & \leq \frac{1}{\pi^2 + 3/4} \| \partial_r f \|_{L^2}^2 + \frac{1}{\pi^2} \| \partial_{x_3} f \|_{L^2}^2 \\
        & \leq \frac{1}{\pi^2} \| \nabla f \|_{L^2}^2.
    \ea \ee
\end{proof}

\subsection{$L^2$ and $L^4$ bounds for $\Gamma$}\label{sec:Gamma}
In this subsection we show that $\|\Gamma\|_{L^2}$ and $\|\Gamma\|_{L^4}$ stay bounded by quantities depending only on the initial data, under the zero-mean condition $\int_D\Gamma_0\,dx=0$.

Before stating the next lemma, we define a cutoff function 
\be
\eta(r):=\sin(\pi r/2), 
\label{eq:eta}
\ee
which is the limit of $\eta_m$ in \cite[Lemma 2.3]{li2026cone}.
Then $\eta(0)=0$, $\eta(1)=1$, $\eta,\partial_r \eta\ge0$, and
\be
0 \leq r\eta \leq 1, \quad 0\le \eta+r \partial_r \eta\le \frac75, \qquad 0\le 2\eta+r \partial_r \eta\le 3\eta+r \partial_r \eta\le \frac{13}{4}.
\label{eq:eta_bound}
\ee

\begin{lemma}
    Suppose $\Gamma$ is a smooth solution to \eqref{eq:GammaU} and \eqref{eq:GammaU_bc} and 
    \be
        \int_D \Gamma_{0} \,dx = 0.
    \ee
    Then for all $t \geq 0$,
    \be
    \| \Gamma \|_{L^2(D)} \leq \|\Gamma_{0} \|_{L^2(D)}.
    \ee
    \label{lemma:L^2_Gamma_U}
\end{lemma}
\begin{proof}
    By testing \eqref{eq:GammaU} with $\Gamma$, we get
    \be
        \frac{1}{2} \ddt \int_D |\Gamma|^2 \, dx + \int_D |\nabla \Gamma |^2 \, dx = \int_0^1 \Gamma^2(1,x_3,t) \, dx_3.
    \ee
    We rewrite the right side as
    \be
        \int_0^1  \Gamma^2(1,x_3,t) \, dx_3 = \int_0^1 \int_0^1 \partial_r( \Gamma^2 r^2 \eta ) \, drdx_3.
    \ee
    Then we can compute
    \be \ba
        \int_0^1 \int_0^1 \partial_r( \Gamma^2 r^2 \eta ) \, drdx_3 &= \int_0^1 \int_0^1 2 \Gamma \partial_r \Gamma \, r^2 \eta  drdx_3 + \int_0^1 \int_0^1 2 \Gamma^2 \, r \eta drdx_3 + \int_0^1 \int_0^1  \Gamma^2 \, r^2 \partial_r \eta drdx_3  \\
        & = \int_0^1\int_0^1 2 \Gamma \partial_r \Gamma r^2 \eta \, drdx_3 + \int_0^1\int_0^1  \Gamma^2r (2\eta + r\partial_r\eta) \, drdx_3 \\
        & \leq \int_0^1\int_0^1 2 |\Gamma \partial_r \Gamma|  \, rdrdx_3 + \frac{13}{4}\int_0^1\int_0^1  \Gamma^2 \, r drdx_3.
        \label{eq:estimate_eta_1}
    \ea \ee
    where in the last inequality sign we have used that $r\eta \leq 1$ and \eqref{eq:eta_bound}. Using Cauchy–Schwarz inequality and Lemma \ref{lemma:Poincare}  we have
    \be \ba
        \int_0^1 \int_0^1 \partial_r( \Gamma^2 r^2 \eta ) \, drdx_3 &\leq \left( \frac{13}{4} + \epsilon \right) \int_0^1 \int_0^1 \Gamma^2 \, rdrdx_3 + \epsilon^{-1} \int_0^1 \int_0^1 |\partial_r \Gamma|^2 \, rdrdx_3 \\
        & \leq \frac{1}{\pi^2} \left( \frac{13}{4} + \epsilon \right) \int_0^1 \int_0^1 |\nabla \Gamma|^2 \, rdrdx_3 + \epsilon^{-1} \int_0^1 \int_0^1 |\nabla \Gamma|^2 \, rdrdx_3.
    \ea \ee
    By setting $\epsilon = \sqrt{19/2}$ we have 
    \be
        \frac{1}{\pi^2} \left( \frac{13}{4} + \sqrt{\frac{19}{2}} \right) + \sqrt{\frac{2}{19}} < 0.97.
    \ee
    Then
    \be 
        \frac{1}{2} \ddt \int_D |\Gamma|^2 \, dx + 0.03 \int_D |\nabla \Gamma |^2 \, dx \leq 0
    \ee
    This implies that for any $t \geq 0$
    \be
    \| \Gamma(t) \|_{L^2(D)} \leq \| \Gamma_{0} \|_{L^2(D)}.
    \ee
\end{proof}

\begin{lemma}
    If $\Gamma$ is a smooth solution to \eqref{eq:GammaU} and \eqref{eq:GammaU_bc} and 
    \be
        \int_D \Gamma_{0} \,dx = 0.
    \ee
    Then for all $t \geq 0$,
    \be
    \|\Gamma\|_{L^4}^4 \leq \max \Big\{    \|\Gamma_{0}\|_{L^4}^4, \frac{4}{3\pi^2} \Big(  4  \frac{87}{5} \Big(\frac{174}{5}\Big)^{3/2} C_S^3 + 87 + \frac{3}{2} \pi^2 \Big) \| \Gamma_{0} \|_{L^2}^4  \Big\}=:G^4.
    \ee
    \label{lemma:L^4_Gamma_U}
\end{lemma}

\begin{proof}
    By testing \eqref{eq:GammaU} with $|\Gamma|^{p-2}\Gamma$, we get
    \be \ba
        & \frac{1}{p} \ddt \| \Gamma\|_{L^p}^p + (p-1) \int_D |\Gamma|^{p-2} |\nabla \Gamma|^2 \, dx 
        \\ & = \frac{1}{p} \ddt \| \Gamma
        \|_{L^p}^p + \frac{4(p-1)}{p^2} \Big\| \nabla |\Gamma|^{p/2} \Big\|_{L^2}^2  = \frac{2(p-1)}{p} \int_0^1 |\Gamma|^p (1,x_3)\,dx_3.
    \ea \ee
    When $p=4$, we get
    \be
        \frac{1}{4} \ddt \|\Gamma\|_{L^4}^4 + \frac{12}{16} \Big\| \nabla \Gamma^{2} \Big\|_{L^2}^2 = \frac{6}{4} \int_0^1 \Gamma^4(1,x_3,t) \, dx_3,
    \ee
    which is the same as
    \be
        \ddt \|\Gamma^2\|_{L^2}^2 + 3 \Big\| \nabla \Gamma^{2} \Big\|_{L^2}^2 = 6\int_0^1 \Gamma^4(1,x_3,t) \, dx_3.
        \label{eq:Gamma_L^4_1}
    \ee
    Using the cut-off function $\eta$ defined in \eqref{eq:eta}, we estimate the right side in the same way as \eqref{eq:estimate_eta_1} to get
    \be \ba
        &\int_0^1 \Gamma^4(1,x_3,t) \, dx_3 = \int_0^1\int_0^1 \partial_r(r^2 \eta (\Gamma^2)^2) \, drdx_3
        \\ & = \int_0^1\int_0^1 2 \Gamma^2 \partial_r \Gamma^2 r^2 \eta \, drdx_3 + \int_0^1\int_0^1 (\Gamma^2)^2 r (2\eta+r\partial_r \eta)\, drdx_3  
        \\ & \leq \frac{1}{4} \int_D |\partial_r (\Gamma^2)|^2 \, dx + 4  \int_D (\Gamma^2)^2 \, dx + \frac{13}{4} \int_D (\Gamma^2)^2 \, dx
        \\ & \leq \frac{1}{4} \|\nabla \Gamma^2 \|_{L^2}^2 + \Big(\frac{13}{4}+4 \Big) \| \Gamma^2\|_{L^2}^2.
    \ea \ee
    Hence we have from \eqref{eq:Gamma_L^4_1} that
    \be
        \ddt \|\Gamma^2\|_{L^2}^2 + \frac{3}{2} \Big\| \nabla \Gamma^{2} \Big\|_{L^2}^2  \leq \frac{87}{2} \| \Gamma^2\|_{L^2}^2.
        \label{eq:Gamma_L^4_2}
    \ee
    Next, let us denote
    \be
        \bar{(\Gamma^2)} := \frac{\|\Gamma^2\|_{L^1}}{|D|} = \frac{\|\Gamma\|_{L^2}^2}{|D|} = 2 \|\Gamma\|_{L^2}^2.
        \label{eq:Gamma_L^4_5}
    \ee
    We then compute 
    \be \ba
        & \| (\Gamma^2) \|_{L^2}^2 =  \| (\Gamma^2) -  \bar{(\Gamma^2)} + \bar{(\Gamma^2)}\|_{L^2}^2 = \int_D  \Big((\Gamma^2) -  \bar{(\Gamma^2)} + \bar{(\Gamma^2)} \Big)^2 \, dx
        \\ & = \| (\Gamma^2) -  \bar{(\Gamma^2)} \|_{L^2}^2 + 2\int_D (\Gamma^2)\bar{(\Gamma^2)}  -  \bar{(\Gamma^2)}^2 \, dx + \| \bar{(\Gamma^2)} \|_{L^2}^2
        \\ & = \| (\Gamma^2) -  \bar{(\Gamma^2)} \|_{L^2}^2 + \| \bar{(\Gamma^2)} \|_{L^2}^2
        \\ & = \| (\Gamma^2) -  \bar{(\Gamma^2)} \|_{L^2}^2 + 2 \| \Gamma \|_{L^2}^4,
        \label{eq:Gamma_L^4_3}
    \ea \ee
    where the cross term vanishes due to 
    \be
        2\int_D (\Gamma^2)\bar{(\Gamma^2)}  -  \bar{(\Gamma^2)}^2 \, dx = 2 \bar{(\Gamma^2)} \int_D (\Gamma^2) - \bar{(\Gamma^2)}\, dx = 2 \bar{(\Gamma^2)} \Big( \int_D \Gamma^2 \, dx - |D|  \bar{(\Gamma^2)} \Big) = 0.
    \ee
    Using H\"older's inequality we have
    \be
        \| (\Gamma^2) -  \bar{(\Gamma^2)} \|_{L^2}^2 \leq \| (\Gamma^2) -  \bar{(\Gamma^2)} \|_{L^1}^{4/5} \| (\Gamma^2) -  \bar{(\Gamma^2)} \|_{L^6}^{6/5}.
        \label{eq:lemma_L^4_Ho}
    \ee
    Since $(\Gamma^2) - \bar{(\Gamma^2)}$ has mean zero, using Sobolev inequality then Lemma \ref{lemma:Poincare} gives,
    \be
        \| (\Gamma^2) -  \bar{(\Gamma^2)} \|_{L^6} \leq C_S \| \nabla (\Gamma^2) \|_{L^2},
        \label{eq:lemma_L^4_Sob}
    \ee
    where $C_S = C_{sob}(1 + 1/\pi)>0$. Then we put \eqref{eq:lemma_L^4_Sob} in \eqref{eq:lemma_L^4_Ho} and use Young's inequality with $\epsilon = 1/58$ to get
    \be \ba
        &  \| (\Gamma^2) -  \bar{(\Gamma^2)} \|_{L^2}^2 \leq  C^{6/5}_S \| \nabla (\Gamma^2) \|_{L^2}^{6/5}\| (\Gamma^2) -  \bar{(\Gamma^2)} \|_{L^1}^{4/5}
        \\ & \leq \frac{1}{58} \| \nabla (\Gamma^2) \|_{L^2}^2 + \frac{2}{5} \Big(\frac{174}{5}\Big)^{3/2} C_S^3 \| (\Gamma^2) -  \bar{(\Gamma^2)} \|_{L^1}^{2}.
    \ea \ee
    Hence from \eqref{eq:Gamma_L^4_2}, Lemma \ref{lemma:L^2_Gamma_U} and the definition \eqref{eq:Gamma_L^4_5}, we have
    \be \ba
         & \ddt \|\Gamma^2\|_{L^2}^2 + \frac{3}{4} \Big\| \nabla \Gamma^{2} \Big\|_{L^2}^2  
         \\ & \leq \frac{87}{2} \frac{2}{5} \Big(\frac{174}{5}\Big)^{3/2} C_S^3 \| (\Gamma^2) -  \bar{(\Gamma^2)} \|_{L^1}^{2} + 87 \|\Gamma \|_{L^2}^4
         \\ & \leq \frac{87}{2} \frac{2}{5} \Big(\frac{174}{5}\Big)^{3/2} C_S^3 \Big( \| (\Gamma^2) \|_{L^1} +  |D| \bar{(\Gamma^2)}  \Big)^{2} + 87 \|\Gamma \|_{L^2}^4
         \\ & \leq 4 \frac{87}{2} \frac{2}{5} \Big(\frac{174}{5}\Big)^{3/2} C_S^3  \| \Gamma \|_{L^2}^4  + 87 \|\Gamma \|_{L^2}^4
         \\ & \leq \Big( 4 \frac{87}{2} \frac{2}{5} \Big(\frac{174}{5}\Big)^{3/2} C_S^3 + 87\Big) \| \Gamma_{0} \|_{L^2}^4.
         \label{eq:Gamma_L^4_4}
    \ea \ee
    Applying Lemma \ref{lemma:Poincare} to $(\Gamma^2) - \bar{(\Gamma^2)}$ and using \eqref{eq:Gamma_L^4_3} we get
    \be
       \| (\Gamma^2) \|_{L^2}^2 - \| \bar{(\Gamma^2)} \|_{L^2}^2 = \| \Gamma^2 - \bar{(\Gamma^2)} \|_{L^2}^2 \leq \frac{1}{\pi^2} \| \nabla \Gamma^2 \|_{L^2}^2.
       \label{eq:Gamma_L^4_6}
    \ee
    Thus \eqref{eq:Gamma_L^4_4} and \eqref{eq:Gamma_L^4_6} implies that
    \be\ba
        & \ddt \|\Gamma^2\|_{L^2}^2 + \frac{3\pi^2}{4} \|\Gamma^2\|_{L^2}^2 \leq \Big( 4 \frac{87}{2} \frac{2}{5} \Big(\frac{174}{5}\Big)^{3/2} C_S^3 + 87 + \frac{3}{2} \pi^2 \Big)\| \Gamma_{0} \|_{L^2}^4.
    \ea \ee
    Applying the Gr\"onwall inequality we get,
    \be
    \|\Gamma^2\|_{L^2}^2 \leq  \|\Gamma_{0}^2\|_{L^2}^2 e^{-\frac{3\pi^2 t}{4}} + \frac{4}{3\pi^2} \Big(  4 \frac{87}{2} \frac{2}{5} \Big(\frac{174}{5}\Big)^{3/2} C_S^3 + 87 + \frac{3}{2} \pi^2 \Big) \| \Gamma_{0} \|_{L^2}^4 \Big( 1 - e^{-\frac{3\pi^2 t}{4}} \Big),
    \ee
    where the right side is a convex combination. Hence we get
    \be
    \|\Gamma^2\|_{L^2}^2 \leq \max \Big\{    \|\Gamma_{0}^2\|_{L^2}^2, \frac{4}{3\pi^2} \Big(  4 \frac{87}{2} \frac{2}{5} \Big(\frac{174}{5}\Big)^{3/2} C_S^3 + 87 + \frac{3}{2} \pi^2 \Big) \| \Gamma_{0} \|_{L^2}^4   \Big\} ,
    \ee
    which is the same as 
    \be
    \|\Gamma\|_{L^4}^4 \leq \max \Big\{    \|\Gamma_{0}\|_{L^4}^4, \frac{4}{3\pi^2} \Big(  4 \frac{87}{2} \frac{2}{5} \Big(\frac{174}{5}\Big)^{3/2} C_S^3 + 87 + \frac{3}{2} \pi^2 \Big) \| \Gamma_{0} \|_{L^2}^4  \Big\}.
    \ee
\end{proof}
Using Lemma \ref{lemma:L^2_Gamma_U} and Lemma \ref{lemma:L^4_Gamma_U}, we have for all $t\ge0$
\be
    \| \Gamma(t) \|_{L^2} \leq \| \Gamma_{0} \|_{L^2}, \qquad \| \Gamma(t) \|_{L^4} \leq G.
    \label{eq:Gamma}
\ee

\subsection{Known estimates for $v_r/r$}\label{sec:known}
We recall some estimates for $v_r/r$ from \cite{chan2026global}, which remain valid under the boundary conditions in \eqref{eq:ASNS}.
\begin{lemma} If $v$ is a smooth solution to \eqref{eq:ASNS}, then
\begin{equation}
    \| \nabla \frac{v_r}{r} \|_{L^2}  \leq \| \Omega \|_{L^2},
    \label{eq:lemma1_estimate1}
\end{equation}
and 
\begin{equation}
    \| \nabla^2 \frac{v_r}{r}\|_{L^2} \leq  \| \partial_{x_3} \Omega \|_{L^2},
    \label{eq:lemma1_estimate2}
\end{equation}
\label{lemma1}
\end{lemma}
\begin{proof}
    The proof can be found in Lemma 2.2 in \cite{chan2026global}. We note that the boundary conditions used in Lemma 2.2 in \cite{chan2026global} are 
    \be
    v_r(1,x_3,t) = \partial_{x_3}v_r(r,1,t) = \partial_{x_3} v_r(r,0,t) = 0,
    \ee 
    and
    $\Omega = 0$ on $\partial D$. These boundary conditions are the same as \eqref{eq:ASNS} and in \eqref{eq:VOmega_bc}.
\end{proof}

\begin{lemma} If $v$ is a smooth solution to \eqref{eq:ASNS}, then for any $r,t$
    \begin{equation}
    \int_0^1 \left( \frac{v_r}{r} \right)(r,x_3,t) \, dx_3 = 0.
    \end{equation}
\label{lemma4}
\end{lemma}
\begin{proof}
The proof can be found in Lemma 2.3 in \cite{chan2026global}, where the boundary conditions used are $\omega_\theta=0$, $\partial_{x_3}v_r = \partial_r v_3 =0 $ on $\partial^H D$ and $v_r = 0$ on $\partial^V D$ which are the same as \eqref{eq:ASNS}.
\end{proof}

\begin{lemma} If $v$ is a smooth solution to \eqref{eq:ASNS} then
\label{lemma:agmon_inequality}
\begin{equation}
    \| \frac{v_r}{r} \|_{L^\infty(D)} \leq C_2 \| \nabla \frac{v_r}{r}\|_{L^2(D)}^{1/2} \| \nabla^2 \frac{v_r}{r} \|_{L^2(D)}^{1/2}
\end{equation}
where $C_2 = (1+\sqrt{2}) 2 \sqrt{\pi} (536)^{1/4}(57452(1+5/\pi^2) + 60)^{1/4} $.
\end{lemma}
\begin{proof}
    The proof can be found in Lemma 2.5 in \cite{chan2026global}. In the proof inside \cite{chan2026global}, an even reflection in $x_3$, which requires $\partial_{x_3}v_r = 0$ on $\partial^H D$, and an odd reflection at $r=1$, which requires $v_r=0$, were used. These conditions are still true under the boundary conditions in \eqref{eq:ASNS}.
\end{proof}

\subsection{A trace inequality for $V$}\label{sec:trace}
In this subsection we estimate the boundary term $\int_0^1 V^4(1,x_3,t)\,dx_3$ that appears in the energy estimate for $V$.

Let 
\be
\phi = \frac{V^2}{r}.
\ee
We defined the weighted average of $\phi$ with the measure $d\mu(r,x_3) = 6r^5 drdx_3$ to be
\be
    \bar{\phi} := \int_D \phi \, d\mu(r,x_3) = \int_0^1\int_0^1 \phi \, 6r^5 drdx_3.
\ee
Due to the definition of $\Gamma = rv_\theta$, $\phi = \Gamma^2/r^4$ and 
\be
\bar{\phi} = 6 \|\Gamma\|_{L^2}^2.
\label{eq:phi_Gamma}
\ee

The next lemma gives a Poincar\'e inequality for $\phi$.
\begin{lemma}
    Let $\phi=\phi(r,x_3)$ be smooth on $(0,1]\times[0,1]$ with $\int_D r^2|\nabla\phi|^2\,dx<\infty$, and let $\bar\phi=\int_0^1\int_0^1 \phi\, 6r^5\,dr\,dx_3$. Then   
    \be \ba
     \int_0^1 |\phi(1,x_3,t) - \bar{\phi}|^2 \, dx_3 \leq  12 \int_D r^2 |\nabla \phi|^2 \, dx.
    \ea \ee
    \label{lemma:mean}
\end{lemma}
\begin{proof}
    For simplicity we suppress the $t$ in the functions. Fix $x_3 \in [0,1]$. Since $\int_0^1\int_0^1  6r^5 drdx_3 = 1$, by Cauchy–Schwarz inequality we have
    \be \ba
        & | \phi(1,x_3) - \bar{\phi} |^2 = \Big| \int_D (\phi(1,x_3) - \phi(r,y)) \, d\mu(r,y) \Big|^2  \leq \int_D |\phi(1,x_3) - \phi(r,y)|^2 \, d\mu(r,y).
    \ea \ee
    We express the right side as
    \be \ba
        & \int_D |\phi(1,x_3) - \phi(r,y)|^2 \, d\mu(r,y) = \int_D |\phi(1,x_3) - \phi(r,x_3) + \phi(r,x_3) - \phi(r,y)|^2 \, d\mu(r,y) 
        \\ & \leq \underbrace{2\int_D |\phi(1,x_3) - \phi(r,x_3)|^2 \, d\mu(r,y)}_{T_1} + \underbrace{2\int_D |\phi(r,x_3) - \phi(r,y)|^2 \, d\mu(r,y)}_{T_2}  
    \ea \ee
    and we estimate $T_1$ and $T_2$. 

    Before estimating $T_1$ we see that
    \be
        |\phi(1,x_3) - \phi(r,x_3)|^2 = \Big( \int_r^1 \partial_r \phi(s,x_3) \,ds 
        \Big)^2
        = \Big( \int_r^1 s^{-3/2} s^{3/2}\partial_r \phi(s,x_3) \,ds 
        \Big)^2
        \leq \frac{1}{2r^2} \int_0^1 s^3 (\partial_r \phi(s,x_3))^2 ds.
    \ee
    Hence
    \be \ba
        \frac{1}{2}T_1 & \leq \int_0^1\int_0^1  \frac{1}{2r^2} \int_0^1 s^3 (\partial_r \phi(s,x_3))^2 \,ds \, 6r^5 drdy
         = 3\int_0^1 r^3 \,dr \int_0^1 s^3 (\partial_r \phi(s,x_3))^2 \, ds 
        \\ & \leq \frac{3}{4} \int_0^1 s^3|\partial_r \phi(s,x_3)|^2 ds.
    \ea \ee
    Similarly, using Cauchy–Schwarz inequality, we have
    \be
        |\phi(r,x_3) - \phi(r,y)|^2 \leq \int_0^1 |\partial_{x_3}\phi(r,z)|^2 \, dz.
    \ee
    Substituting this into $T_2$ gives
    \be \ba
    \frac{1}{2}T_2 \leq \int_0^1\int_0^1 \int_0^1 |\partial_{x_3}\phi(r,z)|^2 \, dz \, 6r^5 drdy \leq 6 \int_D r^2 |\partial_{x_3}\phi|^2 \, dx.
    \ea \ee
    Hence combining the estimate for $T_1$ and $T_2$ we have shown that
    \be
        | \phi(1,x_3) - \bar{\phi} |^2 \leq \frac{6}{4} \int_0^1 s^3|\partial_r \phi(s,x_3)|^2 ds + 12 \int_D r^2 |\partial_{x_3}\phi|^2 \, dx.
    \ee
    Integrating in the $x_3$-direction gives
    \be
       \int_0^1 |\phi(1,x_3) - \bar{\phi}|^2 \, dx_3 \leq  12 \int_D r^2 |\nabla \phi|^2 \, dx.
    \ee
\end{proof} 

Before proving the next lemma, we first do a computation to show that
\be
    \int_D |\nabla V^2|^2 + \frac{V^4}{r^2}\, dx  - \int_0^1 V^4(1,x_3,t)\,dx_3 = \int_D r^2 |\nabla \phi|^2 \, dx.
    \label{eq:W}
\ee
Let $W=V^2$. First note that, $W = V^2 = r\phi$, so $$\partial_r(W)-W/r = r \partial_r \phi$$ and $$\partial_{x_3} W = r \partial_{x_3}\phi.$$ Next let us compute the difference
\be \ba
    & \int_D |\nabla W|^2 + \frac{W^2}{r^2}\, dx  - \Big(\int_D \Big( \partial_r W - \frac{W}{r} \Big)^2 + (\partial_{x_3}W)^2 \, dx \Big)
    \\ &= \int_D \frac{2W \partial_r W}{r} \, dx
    \\ &= \int_0^1 \int_0^1 \partial_r( W^2 ) \,drdx_3
    \\ & = \int_0^1 W^2(1,x_3,t) \, dx_3,
    \label{eq:compute_W_1}
\ea \ee 
where the last line is due to the regularity of the solution $v$: $W=V^2 = v_\theta^2/r = 0$ at $r=0$ \cite[Corollary 1]{liu2009characterization}. Hence \eqref{eq:compute_W_1} implies \eqref{eq:W}

In this paper we will need the estimation of the boundary term $$\int_0^1 V^4(1,x_3,t)\, dx_3.$$ An inequality of the form
\be
\int_0^1 V^4(1,x_3,t) \, dx_3 \leq \gamma \Big( \int_D |\nabla V^2|^2 + \frac{V^4}{r^2} \, dx \Big)
\label{eq:impossible_trace}
\ee
cannot hold for $\gamma<1$. If such inequality is true, then inserting \eqref{eq:W} into \eqref{eq:impossible_trace} gives,
\be \ba
\int_0^1 V^4(1,x_3,t) \, dx_3 \leq \gamma \int_0^1 V^4(1,x_3,t)\, dx_3 + \gamma \int_D r^2 |\nabla \phi|^2 \, dx.
\ea \ee
That is 
\be 
    \int_0^1 V^4(1,x_3,t) \, dx_3 \leq \frac{\gamma}{1-\gamma} \int_D r^2 |\nabla \phi|^2 \, dx.
\ee
Given this inequality we can take $\phi = c$, i.e. $V^2 = cr$, where $c>0$, which correspond to $v_\theta= \pm \sqrt c\,r$ and it satisfies the total--slip--boundary condition. Then we get
\be
c^2 \leq 0,
\ee
which is impossible. 

On the other hand, \eqref{eq:W} shows that \eqref{eq:impossible_trace} holds with $\gamma=1$, with equality if and only if $\nabla\phi \equiv 0$. Since $D$ is connected, this means $\phi \equiv c$ for some constant $c\ge0$, i.e. $V^2 = cr$ and $v_\theta = \pm\sqrt{c}\,r$. Taking $c>0$ shows that $\gamma = 1$ is sharp.

However, we can split the boundary value in terms of its mean and the fluctuations, that is
\be
\phi(1,x_3,t) = \bar{\phi} + (\phi(1,x_3,t) - \bar{\phi})
\ee
where the fluctuations is controlled by the Poincar\'e inequality in Lemma \ref{lemma:mean} and the relation \eqref{eq:W}. The mean is equal to $6\|\Gamma\|_{L^2}^2$ as shown in \eqref{eq:phi_Gamma}. This motivates the next lemma.
\begin{lemma}
\label{lemma:trace}
    We have
    \be
    \int_0^1 V^4(1,x_3,t) dx_3 \leq \frac{24}{25} \Big( \int_D |\nabla V^2|^2 + \frac{V^4}{r^2} \, dx \Big) + \frac{72}{25} \| \Gamma \|_{L^2}^4.
    \label{eq:trace}
    \ee
\end{lemma}
\begin{proof}
    We first start with the boundary term and using Lemma \ref{lemma:mean}: 
    \be \ba
        \int_0^1 V^4(1,x_3,t)\,dx_3 & = \int_0^1 \phi^2(1,x_3,t) \, dx_3 
        \\ & = \int_0^1 (\phi(1,x_3,t) - \bar{\phi} + \bar{\phi})^2 \, dx_3
        \\ & \leq 2\int_0^1 |\phi(1,x_3,t)-\bar{\phi}|^2 \, dx_3 + 2 \int_0^1 |\bar{\phi}|^2 \, dx_3 
        \\ & \leq 24 \int_D r^2 |\nabla \phi|^2 \, dx + 2|\bar{\phi}|^2.
    \ea \ee
    By \eqref{eq:W}, $\int_D r^2|\nabla\phi|^2\,dx = \int_D |\nabla V^2|^2 + \frac{V^4}{r^2}\,dx - \int_0^1 V^4(1,x_3,t)\,dx_3$. Substituting this and moving the boundary term to the left side, we have
    \be \ba
        25 \int_0^1 V^4(1,x_3,t) \, dx_3 \leq 24 \int_D |\nabla V^2|^2 + \frac{V^4}{r^2}\, dx + 2|\bar{\phi}|^2.
    \ea \ee
    Recall from \eqref{eq:phi_Gamma} that $\bar{\phi}=6\|\Gamma\|_{L^2}^2$. We have
    \be \ba
        25 \int_0^1 V^4(1,x_3,t) \, dx_3 \leq 24 \int_D |\nabla V^2|^2 + \frac{V^4}{r^2}\, dx + 72\|\Gamma\|_{L^2}^4.
    \ea \ee
    This gives \eqref{eq:trace}.
\end{proof}

\section{Proof of Theorem \ref{thm:main}} \label{sec:proof}
In this section we give a proof of Theorem \ref{thm:main}. We first start with the energy estimate of $V$ and $\Omega$. Then we estimate the $L^4$ norm of $b = v_re_r + v_3 e_3$ and the $L^4$ norm of $v_\theta$. This will be enough to guarantee the regularity of the solution by the Ladyzhenskaya--Prodi--Serrin regularity criterion.

\subsection{Energy estimate for $V$ and $\Omega$}
Since $v_0$ is smooth, the local existence of a strong solution to $\eqref{eq:ASNS}$ under the conditions of Theorem \ref{thm:main} follows from the classical theory. For example one can follow the proof in \cite{chen2010navier} by reflecting to the smooth periodic cylinder, and taking $1/\zeta = 0$ in the paper, where $\zeta$ is the slip length in \cite{chen2010navier}. So we assume there exists a strong solution to $\eqref{eq:ASNS}$ in $D \times [0,T^*)$ with initial data $v_0$, where $T^*$ be the maximal time of existence of this strong solution. We will derive bounds on $[0,T^*)$ that do not depend on $T^*$; the strong solution can then be continued for all time.

We multiply the first equation of \eqref{eq:VOmega} by $V^3$ and integrate in space. The transport term vanishes since $\operatorname{div} b = 0$, $v_r = 0$ on $\partial^V D$ and $v_3=0$ on $\partial^H D$, and we get
\be \ba
 &\frac{1}{4} \ddt \|V^2 \|_{L^2}^2 + \frac{3}{4} \int_D |\nabla V^2|^2 \, dx  + \frac{3}{4} \int_D \frac{1}{r^2} V^4 \, dx \\ &= \int_D \nabla \cdot(V^3 \nabla V) \, dx - \frac{3}{2} \int_D \frac{v_r}{r} V^4 \, dx + \frac{1}{4} \int_0^1 V^4(1,x_3,t) \, dx_3 ,
\ea \ee
where the last term comes from $\int_D \frac1r \partial_r V\, V^3\,dx = \frac14\int_0^1\int_0^1 \partial_r(V^4)\,drdx_3$ and $V = O(r^{1/2})$ near the axis by \cite[Corollary 1]{liu2009characterization}. By the boundary conditions \eqref{eq:VOmega_bc}, we have
\be
    \int_D \nabla \cdot(V^3\nabla V) \, dx = \int_0^1 V^3 \partial_r V (1,x_3,t)\, dx_3 = \frac{1}{2} \int_0^1 V^4(1,x_3,t) \, dx_3,
\ee
hence,
\be \ba
\label{eq:V_es_1}
 &\frac{1}{4} \ddt \|V^2 \|_{L^2}^2 + \frac{3}{4} \int_D |\nabla V^2|^2 \, dx  + \frac{3}{4} \int_D \frac{1}{r^2} V^4 \, dx =  \frac{3}{4} \int_0^1 V^4(1,x_3,t) \, dx_3 - \frac{3}{2} \int_D \frac{v_r}{r} V^4 \, dx .
\ea \ee
By Lemma \ref{lemma:trace} and \eqref{eq:Gamma} we have
\be \ba
    &\frac{1}{4} \ddt \|V^2 \|_{L^2}^2 + \frac{3}{100} \int_D |\nabla V^2|^2 \, dx  + \frac{3}{100} \int_D \frac{1}{r^2} V^4 \, dx
  \leq \frac{54}{25} \|\Gamma_0\|_{L^2}^4 \underbrace{- \frac{3}{2} \int_D \frac{v_r}{r} V^4 \, dx}_{T_1}.
 \label{eq:V_es_2}
\ea \ee

Next, we estimate $T_1$. Since $V^4 = \big(\frac{V^2}{r}\big)^{3/2} |\Gamma|$, using H\"older's inequality, Lemma \ref{lemma1}, Lemma \ref{lemma:agmon_inequality}, Young's inequality and \eqref{eq:Gamma}, we have
\be \ba
T_1 & \leq  \frac{3}{2} \int_D \frac{|v_r|}{r} \Big( \frac{V^2}{r} \Big)^{3/2} |\Gamma| \,dx \leq \frac{3}{2} \Big\|\frac{v_r}{r} \Big\|_{L^\infty}  \Big\|\frac{V^2}{r} \Big\|_{L^2}^{3/2} \|\Gamma\|_{L^4}
\\ & \leq C_2^4 \frac{2187}{4096} \frac{1}{\delta_1^3} \|\Omega\|_{L^2}^2 \|\partial_{x_3} \Omega \|_{L^2}^2 \|\Gamma\|_{L^4}^4 + \delta_1 \Big\| \frac{V^2}{r} \Big\|_{L^2}^2
\\ & \leq C_2^4 \frac{2187}{4096} \frac{1}{\delta_1^3} \|\Omega\|_{L^2}^2 \|\partial_{x_3} \Omega \|_{L^2}^2\, G^4 + \delta_1 \Big\| \frac{V^2}{r} \Big\|_{L^2}^2 .
\label{eq:T1}
\ea \ee
Combining \eqref{eq:V_es_2} and \eqref{eq:T1} we have
\be \ba
    &\frac{1}{4} \ddt \|V^2 \|_{L^2}^2 + \frac{3}{100} \int_D |\nabla V^2|^2 \, dx  + \Big(\frac{3}{100} - \delta_1\Big) \int_D \frac{1}{r^2} V^4 \, dx
 \\ & \leq \frac{54}{25} \|\Gamma_0\|_{L^2}^4 + C_2^4\frac{2187}{4096} \frac{1}{\delta_1^3} \|\Omega\|_{L^2}^2 \|\partial_{x_3} \Omega \|_{L^2}^2\, G^4 .
 \label{eq:V_es_3}
\ea \ee

Using $\Omega$ as a test function in the second equation of \eqref{eq:VOmega}, we obtain
\be \ba
\frac{1}{2} \ddt \| \Omega \|_{L^2}^2 + \| \nabla \Omega \|_{L^2}^2 + \int_0^1 \Omega^2(0,x_3,t)\,dx_3 &= \int_D \frac{\partial_{x_3}V^2}{r} \Omega \, dx = - \int_D  \frac{V^2}{r} \partial_{x_3}\Omega \, dx
\\ & \leq \delta_3 \Big\|\frac{V^2}{r}\Big\|_{L^2}^2 + \frac{1}{4\delta_3} \|\partial_{x_3}\Omega\|_{L^2}^2,
\label{eq:Omega_es_1}
\ea \ee
where we integrated by parts in $x_3$ and used $\Omega = 0$ on $\partial^H D$.

Now we pick
\be
\delta_1 = \frac{1}{400}, \qquad \delta_3 = \frac{1}{3}.
\ee
Adding $\frac{1}{20}$\eqref{eq:Omega_es_1} and \eqref{eq:V_es_3}, dropping the nonnegative term $\int_0^1 \Omega^2(0,x_3,t)\,dx_3$ and using $\|\partial_{x_3}\Omega\|_{L^2} \le \|\nabla \Omega\|_{L^2}$, we get
\be \ba
    &\frac{1}{4} \ddt \|V^2 \|_{L^2}^2 + \frac{1}{40} \ddt \| \Omega \|_{L^2}^2+ \frac{3}{100} \int_D |\nabla V^2|^2 \, dx + \frac{1}{80} \| \nabla \Omega \|_{L^2}^2 + \frac{13}{1200} \int_D \frac{1}{r^2} V^4 \, dx
 \\ & \leq  \frac{54}{25} \| \Gamma_{0} \|_{L^2}^4 + \frac{2187}{4096}\, 400^3\, C_2^4\, G^4\, \|\Omega\|_{L^2}^2 \| \nabla \Omega \|_{L^2}^2 ,
 \label{eq:ES_0}
\ea \ee
where $\frac{3}{100} - \frac{1}{400} - \frac{1}{20}\cdot\frac13 = \frac{13}{1200}$ and $\frac{1}{20}\big(1 - \frac34\big) = \frac{1}{80}$.

We now close the estimate by a continuity argument. Let $f(t) = \frac{1}{4} \|V^2(t)\|_{L^2}^2 + \frac{1}{40} \|\Omega(t)\|_{L^2}^2$ and
\be
T^{**} := \sup\Big\{ t \in [0,T^*) : \|\Omega(s)\|_{L^2}^2 \le 80 F \text{ for all } s \in [0,t] \Big\}.
\ee
Since $\|\Omega(0)\|_{L^2}^2 \le 40 f(0) \le 40F < 80 F$, continuity gives $T^{**} > 0$. For $t \in [0,T^{**})$, condition (3) of Theorem \ref{thm:main} gives
\be
\frac{2187}{4096}\, 400^3\, C_2^4\, G^4\, \|\Omega\|_{L^2}^2 \le 80\cdot \frac{2187}{4096}\, 400^3\, C_2^4\, F\, G^4 \le \frac{1}{160},
\ee
so the last term in \eqref{eq:ES_0} is absorbed and
\be \ba
 & \frac{1}{4} \ddt \|V^2 \|_{L^2}^2 + \frac{1}{40} \ddt \| \Omega \|_{L^2}^2+ \frac{3}{100} \int_D |\nabla V^2|^2 \, dx + \frac{1}{160} \| \nabla \Omega \|_{L^2}^2 + \frac{13}{1200} \int_D \frac{1}{r^2} V^4 \, dx
  \leq \frac{54}{25} \| \Gamma_{0} \|_{L^2}^4 .
 \label{eq:ES_1}
\ea \ee

Next we note that $\Omega$ satisfies the Poincar\'e inequality $\|\Omega\|_{L^2} \leq \|\nabla\Omega\|_{L^2}$. Indeed, using $\Omega(1,x_3,t)=0$,
\be
\int_0^1\Omega^2r\,dr=-\int_0^1\Omega\,\partial_r\Omega\,r^2\,dr\le\Big(\int_0^1\Omega^2r\,dr\Big)^{1/2}\Big(\int_0^1(\partial_r\Omega)^2r^3\,dr\Big)^{1/2},
\ee
so $\int_0^1 \Omega^2 r\,dr \le \int_0^1 (\partial_r\Omega)^2 r^3\,dr \le \int_0^1 (\partial_r \Omega)^2 r\,dr$ since $r \le 1$; integrating in $x_3$ gives the Poincar\'e inequality.

Since $r\leq 1$, we have $\frac{13}{1200}\| r^{-1}V^2 \|_{L^2}^2 \ge \frac{1}{100}\| V^2 \|_{L^2}^2 = \frac{1}{25}\cdot \frac14 \|V^2\|_{L^2}^2$, and by the Poincar\'e inequality $\frac{1}{160}\|\nabla\Omega\|_{L^2}^2 \ge \frac{1}{160}\|\Omega\|_{L^2}^2 \ge \frac{1}{25}\cdot\frac{1}{40}\|\Omega\|_{L^2}^2$. Hence \eqref{eq:ES_1} implies
\be
f' + \frac{1}{25} f \leq \frac{54}{25} \| \Gamma_{0} \|_{L^2}^4 \qquad \text{on } [0,T^{**}).
\ee
By the Gr\"onwall inequality,
\be
 f(t) \leq  f(0)e^{-\frac{t}{25}} + 54 \| \Gamma_{0} \|_{L^2}^4 \big( 1- e^{-\frac{t}{25}} \big).
\label{eq:Gronwall}
\ee
The right side is a convex combination of $f(0)$ and $54 \| \Gamma_{0} \|_{L^2}^4$, hence
\be
 \frac{1}{4} \|V^2(t)\|_{L^2}^2 + \frac{1}{40} \|\Omega(t)\|_{L^2}^2 = f(t) \leq \max \Big\{ f(0),\ 54 \| \Gamma_{0} \|_{L^2}^4 \Big\} = F, \qquad t \in [0,T^{**}).
 \label{eq:f_bound}
\ee
In particular $\|\Omega(t)\|_{L^2}^2 \le 40 F < 80 F$ on $[0,T^{**})$. If $T^{**} < T^*$, then by continuity $\|\Omega(t)\|_{L^2}^2 < 80F$ on a neighborhood of $T^{**}$, contradicting the definition of $T^{**}$. Therefore $T^{**} = T^*$, and \eqref{eq:f_bound} holds on $[0,T^*)$. We emphasize that the bound $F$ does not depend on $T^*$.

\subsection{$L^4$ bound for $b$}
Using the formula for the gradient of $b$,
\be
|\nabla b|^2=(\partial_rv_r)^2+(\partial_{x_3}v_r)^2+(\partial_rv_3)^2+(\partial_{x_3}v_3)^2+\frac{v_r^2}{r^2},
\ee
and
\be
\omega_\theta=\partial_{x_3}v_r-\partial_rv_3,\qquad 0=\operatorname{div}b=\partial_rv_r+\frac{v_r}{r}+\partial_{x_3}v_3 ,
\ee
we can compute
\be
\int_D|\nabla b|^2 \, dx -\int_D\omega_\theta^2 \, dx = \int_D |\operatorname{div}b|^2\,dx-2\int_D\Big(\partial_rv_r\frac{v_r}{r}+\partial_rv_r\,\partial_{x_3}v_3+\frac{v_r}{r}\partial_{x_3}v_3-\partial_{x_3}v_r\,\partial_rv_3\Big)dx .
\ee
The first integral vanishes since $\operatorname{div} b = 0$. We will show that each remaining term reduces to a boundary term, and all of them vanish.

First,
\be
\int_D2\partial_rv_r\frac{v_r}{r}\,dx=\int_0^1\!\!\int_0^1\partial_r(v_r^2)\,dr\,dx_3=\int_0^1 v_r^2\Big|_{r=0}^{r=1}dx_3=0,
\ee
because $v_r=0$ at $r=1$ and $v_r = O(r)$ near the axis $r=0$ \cite[Corollary 1]{liu2009characterization}.

Next, we notice that
\be
\partial_rv_r\,\partial_{x_3}v_3-\partial_{x_3}v_r\,\partial_rv_3=\partial_r(v_r\partial_{x_3}v_3)-\partial_{x_3}(v_r\partial_rv_3).
\ee
We compute
\be \ba
& \int_D\Big(\partial_rv_r\partial_{x_3}v_3-\partial_{x_3}v_r\partial_rv_3+\frac{v_r}{r}\partial_{x_3}v_3\Big) \, dx
\\ & = \int_0^1\!\!\int_0^1 \Big( r\partial_r(v_r\partial_{x_3}v_3) - r\partial_{x_3}(v_r\partial_rv_3) \Big) \, dr\,dx_3 + \int_0^1\!\! \int_0^1 v_r\partial_{x_3}v_3 \, dr\,dx_3
\\ & = \int_0^1\!\!\int_0^1 \Big( \partial_r( r v_r \partial_{x_3}v_3) - v_r  \partial_{x_3}v_3 - r\partial_{x_3}(v_r\partial_rv_3) \Big) \, dr\,dx_3 + \int_0^1\!\! \int_0^1 v_r\partial_{x_3}v_3 \, dr\,dx_3
\\ & = \int_0^1\!\!\int_0^1 \Big( \partial_r( r v_r \partial_{x_3}v_3) - \partial_{x_3}(r v_r\partial_rv_3) \Big) \, dr\,dx_3
\\ & =\int_0^1 \big[r\,v_r\partial_{x_3}v_3\big]_{r=0}^{r=1}\, dx_3-\int_0^1 \big[r\,v_r\partial_rv_3\big]_{x_3=0}^{x_3=1}\, dr = 0,
\ea \ee
where all the boundary terms vanish because $v_r=0$ at $r=1$, the factor $r$ vanishes at the axis, and $v_3=0$ on $\partial^H D$, so $\partial_r v_3=0$ on $\partial^H D$.

Hence, since $\omega_\theta=r\Omega$ and $r\le1$,
\be
\|\nabla b\|_{L^2}=\|\omega_\theta\|_{L^2}\le\|\Omega\|_{L^2}.
\ee
Moreover, since $r \le 1$, $\|v_r\|_{L^2}\le\|v_r/r\|_{L^2}$. For $v_3$, fix $r$ and write
\be
\int_0^1 v_3^2 \, dx_3 = \int_0^1 v_3^2 \, \partial_{x_3} \Big(x_3 - \frac{1}{2} \Big) \, dx_3.
\ee
Since $v_3=0$ on $\partial^H D$, the boundary term vanishes after integrating by parts, and using $|x_3 - \frac12| \le \frac12$ and the Cauchy--Schwarz inequality gives
\be
\int_0^1 v_3^2 \, dx_3 = - 2\int_0^1 v_3 \partial_{x_3}v_3 \Big(x_3 - \frac{1}{2} \Big)\, dx_3 \leq \Big( \int_0^1 v_3^2 \, dx_3 \Big)^{1/2}  \Big( \int_0^1 |\partial_{x_3}v_3|^2 \, dx_3 \Big)^{1/2},
\ee
which means
\be
     \int_0^1 v_3^2 \, dx_3  \leq  \int_0^1 |\partial_{x_3}v_3|^2 \, dx_3 .
\ee
Multiplying by $r$ and integrating in $r$ gives $\|v_3\|_{L^2}\le\|\partial_{x_3} v_3\|_{L^2}$. Since $v_r/r$ and $\partial_{x_3}v_3$ are distinct terms in $|\nabla b|^2$,
\be
\|b\|_{L^2}^2 = \|v_r\|_{L^2}^2 + \|v_3\|_{L^2}^2 \le \Big\|\frac{v_r}{r}\Big\|_{L^2}^2 + \|\partial_{x_3}v_3\|_{L^2}^2 \le \|\nabla b\|_{L^2}^2 .
\ee
Therefore $\|b\|_{H^1}^2 \le 2\|\nabla b\|_{L^2}^2 \le 2\|\Omega\|_{L^2}^2 \le 80 F$ on $[0,T^*)$ by \eqref{eq:f_bound}. Applying the Sobolev inequality and H\"older's inequality, we obtain
\be
\sup_{t\in[0,T^*)} \|b(t)\|_{L^4} \le |D|^{1/12} C_{sob} \sqrt{80F}.
\ee

\section{Completion of the proof of Theorem \ref{thm:main}} \label{sec:completion}
By \eqref{eq:f_bound} and $r \le 1$,
\be
\|v_\theta\|_{L^4}^4 = \int_D v_\theta^4\,dx \le \int_D \frac{v_\theta^4}{r^2}\,dx = \|V^2\|_{L^2}^2 \le 4F .
\ee
Together with the bound on $b$, this shows that $v_r$, $v_\theta$ and $v_3$ are all in $L^\infty(0,T^*;L^4(D))$, with a bound depending only on $F$ and not on $T^*$. By the continuation principle the solution exists for all time, and $v \in L^\infty(0,\infty;L^4(D))$. Hence by the Ladyzhenskaya--Prodi--Serrin regularity criterion \cite{serrin1962interior,prodi1959teorema,ladyzhenskaya1967uniqueness}, $v$ is regular up to any finite time $T$ in the interior of $D$.

Also, a blow-up can only occur on $\{r=0\}$ by the result in \cite{caffarelli1982partial}, so it remains to rule out blow-up on $\{ r=0 \} \cap \partial^H D$. The total-slip condition is the generalized Navier condition of \cite{neustupa2018regularity} with zero friction, and by the result in \cite{neustupa2018regularity} the interior regularity criterion also implies regularity on $\partial^H D$. So $v$ is regular up to any finite time $T$ in $D$.  \qed

\section*{Acknowledgments} We are deeply grateful to Professors Qi S. Zhang, Zijin Li, Xinghong Pan, Xin Yang, Na Zhao and Daoguo Zhou for the helpful discussions and suggestions.

\bibliographystyle{plain}

\end{document}